\documentclass[12pt,twoside]{amsart}
\usepackage{amscd}      
\usepackage{amssymb}
\usepackage{amsmath}
\usepackage{xypic}      
\LaTeXdiagrams          
\usepackage[all,v2]{xy}
\xyoption{2cell} \UseAllTwocells \xyoption{frame} \CompileMatrices
\allowdisplaybreaks[3]

\usepackage{latexsym}
\usepackage{epsfig}
\usepackage{amsfonts}
\usepackage{enumerate}

\newtheorem{prop}{Proposition}[section]
\newtheorem{lem}[prop]{Lemma}

\newtheorem{cor}[prop]{Corollary}
\newtheorem{thm}[prop]{Theorem}

\numberwithin{equation}{section}

\newtheorem{defn}[prop]{Definition}

\newtheorem{rmk}[prop]{Remark}

\newenvironment{pf}{\begin{trivlist}\item[]{\sc Proof.}}%
            {\nolinebreak $\Box$ \end{trivlist}}

\newcommand{\noprint}[1]{}

\newcommand{\ldiag}[1]%
       {\makebox[0cm]{${\scriptstyle#1}\downarrow\phantom{\scriptstyle#1}$}}
\newcommand{\ldiagup}[1]%
       {\makebox[0cm]{${\scriptstyle#1}\uparrow\phantom{\scriptstyle#1}$}}
\newcommand{\rdiag}[1]%
       {\makebox[0cm]{$\phantom{\scriptstyle#1}\downarrow{\scriptstyle#1}$}}
\newcommand{\sediagr}[1]%
       {\makebox[0cm]{$\phantom{\scriptstyle#1}\searrow{\scriptstyle#1}$}}
\newcommand{\nediagr}[1]%
       {\makebox[0cm]{$\phantom{\scriptstyle#1}\nearrow{\scriptstyle#1}$}}
\newcommand{\rdiagup}[1]%
       {\makebox[0cm]{$\phantom{\scriptstyle#1}\uparrow{\scriptstyle#1}$}}
\newcommand{\swdiag}[1]%
       {\makebox[0cm]{$\phantom{\scriptstyle#1}\swarrow{\scriptstyle#1}$}}
\newcommand{\sediag}[1]%
       {\makebox[0cm]{${\scriptstyle#1}\searrow\phantom{\scriptstyle#1}$}}
\newcommand{\nediag}[1]%
       {\makebox[0cm]{${\scriptstyle#1}\nearrow\phantom{\scriptstyle#1}$}}

\newcommand{\doublearrowstack}[2]%
                      {{{{\scriptstyle#1}\atop{\textstyle\longrightarrow}}\atop{{\textstyle\longrightarrow}\atop{\scriptstyle#2}}}}
\newcommand{\rightleftarrowstack}[2]%
                      {{{{\scriptstyle#1}\atop{\textstyle\longrightarrow}}\atop{{\textstyle\longleftarrow}\atop{\scriptstyle#2}}}}
\newcommand{\leftrightarrowstack}[2]%
                      {{{{\scriptstyle#1}\atop{\textstyle\longleftarrow}}\atop{{\textstyle\longrightarrow}\atop{\scriptstyle#2}}}}

\newcommand{\overtoparrow}%
{\makebox[0cm]{\beginpicture \setcoordinatesystem units
<.8cm,.4cm> point at 0 0 \setplotarea x from -3 to 3, y from 0 to
1 \setquadratic \plot -3 0 0 1 3 0 / \put{\vector(3,-1){0}}[Bl] at
3 0
\endpicture}}

\newcommand{\underbottomarrow}%
{\makebox[0cm]{\beginpicture \setcoordinatesystem units
<.8cm,.4cm> point at 0 0 \setplotarea x from -3 to 3, y from 0 to
1 \setquadratic \plot -3 1 0 0 3 1 / \put{\vector(3,1){0}}[Bl] at
3 1
\endpicture}}

\newcommand{\ses}[5]%
{0\longrightarrow#1\stackrel{#2}{ \longrightarrow}#3\stackrel{#4}{
\longrightarrow}#5\longrightarrow0}

\newcommand{\dt}[6]%
{#1\stackrel{#2}{longrightarrow}#3
\stackrel{#4}{\longrightarrow}#5 \stackrel{#6}{\longrightarrow}
#1[1]}

\newcommand{\cat}[1]%
{(\mbox{\rm #1})}

\def\Label#1{\label{#1}{\tt [#1]}\phantom{h}}
\def\Label{\label}

\title[Semi-projective toric DM stacks]{Note on Orbifold Chow Ring of Semi-projective Toric Deligne-Mumford Stacks}

\author{Yunfeng Jiang}
\address{Department of Mathematics\\ University of British Columbia\\ 1984 Mathematics Road\\Vancouver\\ BC V6T 1Z2\\ Canada}
\email{jiangyf@math.ubc.ca}

\author{Hsian-Hua Tseng}
\address{Department of Mathematics\\ University of British Columbia\\ 1984 Mathematics Road\\Vancouver\\ BC V6T 1Z2\\ Canada}
\email{hhtseng@math.ubc.ca}

\date{\today}

\begin{document}
\begin{abstract}
We prove a formula for the orbifold Chow ring of semi-projective toric DM stacks, generalizing the orbifold Chow ring formula of projective toric DM stacks by Borisov-Chen-Smith. We also consider a special kind of semi-projective toric DM stacks, the Lawrence toric DM stacks. We prove that the orbifold Chow ring of a Lawrence toric DM stack is isomorphic to the orbifold Chow ring of its associated hypertoric DM stack studied in \cite{JT}.
\end{abstract}

\maketitle

\section{Introduction}

The main goal of this paper is to generalize the orbifold Chow ring formula of Borisov-Chen-Smith for 
projective toric DM stacks to the case of semi-projective toric DM stacks.

In the paper \cite {BCS}, Borisov, Chen, and Smith developed the theory of toric DM stacks using stacky fans. A stacky fan is a triple $\mathbf{\Sigma}=(N,\Sigma,\beta)$, where $N$ is a finitely generated abelian group, $\Sigma$ is a simplicial fan in the lattice $\overline{N}:=N\slash \text{torsion}$ and $\beta: \mathbb{Z}^n\rightarrow N$ is a map given by a collection of vectors $\{b_1,\cdots,b_n\}\subset N$ such that the images $\{\overline{b}_{1},\cdots,\overline{b}_{n}\}$ generate the fan $\Sigma$.  A toric DM stack $\mathcal{X}(\mathbf{\Sigma})$ is defined using $\mathbf{\Sigma}$; it is a quotient stack whose coarse moduli space is the toric variety $X(\Sigma)$ corresponding to the simplicial fan $\Sigma$.

The construction of toric DM stacks was slightly generalized later in \cite{Jiang}, in which the notion of extended stacky fans was introduced. This new notion is based on that of stacky fans plus some extra data. Extended stacky fans yield toric DM stacks in the same way as stacky fans do. The main point is that extended stacky fans provide presentations of toric DM stacks not available from stacky fans.

When $X(\Sigma)$ is projective, it is found in \cite{BCS} that the orbifold Chow ring (or Chen-Ruan cohomology ring) of $\mathcal{X}(\mathbf{\Sigma})$ is isomorphic to a deformed ring of the group ring of $N$. We call a toric DM stack $\mathcal{X}(\mathbf{\Sigma})$ {\em semi-projective} if its coarse moduli space $X(\Sigma)$ is semi-projective. Hausel and Sturmfels \cite{HS} computed the Chow ring of semi-projective toric varieties. Their answer is also known as the ``Stanley--Reisner'' ring of a fan. Using their result, we prove a formula of the orbifold Chow ring of semi-projective toric DM stacks.

Consider an extended stacky fan $\mathbf{\Sigma}=(N,\Sigma,\beta)$, where $\Sigma$ is the simplicial fan of the semi-projective  toric variety $X(\Sigma)$. Let $N_{tor}$ be the torsion subgroup of $N$, then $N=\overline{N}\oplus N_{tor}$. Let  $N_{\Sigma}:=|\Sigma|\oplus N_{tor}$. Note that $|\Sigma|$ is convex, so $|\Sigma|\oplus N_{tor}$ is a subgroup of $N$. Define the deformed ring $\mathbb{Q}[N_{\Sigma}]:=\bigoplus_{c\in N_{\Sigma}}\mathbb{Q}y^{c}$ with the product structure given by
\begin{equation}\Label{productA}
y^{c_{1}}\cdot y^{c_{2}}:=\begin{cases}y^{c_{1}+c_{2}}&\text{if
there is a cone}~ \sigma\in\Sigma ~\text{such that}~ \overline{c}_{1}\in\sigma, \overline{c}_{2}\in\sigma\,;\\
0&\text{otherwise}\,.\end{cases}
\end{equation}
Note that if $\mathcal{X}(\mathbf{\Sigma})$ is projective, then $N_{\Sigma}=N$ and $\mathbb{Q}[N_{\Sigma}]$ is the deformed ring
$\mathbb{Q}[N]^{\mathbf{\Sigma}}$ in \cite{BCS}.
Let $A^{*}_{orb}(\mathcal{X}(\mathbf{\Sigma}))$ denote the orbifold Chow ring of the toric DM stack $\mathcal{X}(\mathbf{\Sigma})$.

\begin{thm}\label{main}
Assume that $\mathcal{X}(\mathbf{\Sigma})$ is semi-projective. There is an isomorphism of rings
$$A^{*}_{orb}(\mathcal{X}(\mathbf{\Sigma}))\cong \frac{\mathbb{Q}[N_{\Sigma}]}{\{\sum_{i=1}^{n}e(b_{i})y^{b_{i}}:e\in N^{\star}\}}.$$
\end{thm}

The strategy of proving Theorem \ref{main} is as follows. We use a formula in \cite{HS} for the ordinary Chow ring of semi-projective toric varieties. We prove that each twisted sector is also a semi-projective toric DM stack. With this, we use a method similar to that in \cite{BCS} and \cite{Jiang} to prove the isomorphism as modules. The argument to show the isomorphism as rings is the same as that in \cite{BCS}, except that we only take elements in the support of the fan. 

An interesting class of examples of semi-projective toric DM stack is the Lawrence toric DM stacks. We discuss the properties of such stacks. We prove that each 3-twisted sector or twisted sector is again a Lawrence toric DM stack. This allows us to draw connections to hypertoric DM stacks studied in \cite{JT}. We prove that the orbifold Chow ring of a Lawrence toric DM stack is isomorphic to the orbifold Chow ring of its associated hypertoric DM stack. This is an analog of Theorem 1.1 in \cite{HS} for orbifold Chow rings.

The rest of this text is organized as follows. In Section \ref{semi-pro} we define semi-projective toric DM stacks and prove Theorem 1.1. Results on Lawrence toric DM stacks are discussed in Section \ref{Lawrence}.
\subsection*{Conventions}
In this paper we work entirely algebraically over the field of
complex numbers. Chow rings and orbifold Chow rings are taken with
rational coefficients. By an orbifold we mean a smooth
Deligne-Mumford stack with trivial generic stabilizer.

For a simplicial fan $\Sigma$, we use $|\Sigma|$ to represent the lattice points in $\Sigma$.
Note that if $\Sigma$ is convex, $|\Sigma|$ is a free abelian subgroup of $N$.
We write $N^{\star}$ for $Hom_{\mathbb{Z}}(N,\mathbb{Z})$ and $N\to \overline{N}$ the natural map of modding out torsions. We refer to \cite{BCS} for the construction of the Gale dual
$\beta^{\vee}: \mathbb{Z}^{m}\to DG(\beta)$ of $\beta: \mathbb{Z}^{m}\to N$. 

\subsection*{Acknowledgments}
We thank Kai Behrend and Nicholas Proudfoot for valuable discussions.
\section{Semi-projective toric DM stacks and their orbifold Chow rings}\label{semi-pro}
In this section we define semi-projective toric DM stacks and discuss their properties.

\subsection{Semi-projective toric DM stacks}

\begin{defn}[\cite{HS}]
A  toric variety $X$ is called semi-projevtive if the natural map 
$$\pi: X\rightarrow X_{0}=\text{Spec}(H^{0}(X,\mathcal{O}_{X})),$$
is projective and $X$ has at least one torus-fixed point.
\end{defn}

\begin{defn}[\cite{Jiang}]\label{stackyfan}
An extended stacky fan $\mathbf{\Sigma}$ is a triple $(N,\Sigma,\beta)$,
where  $N$ is a finitely generated abelian group, $\Sigma$ is a
simplicial fan  in $N_{\mathbb{R}}$ and $\beta:
\mathbb{Z}^{m}\to N$ is the map determined by the
elements $\{b_{1},\cdots,b_{m}\}$ in $N$ such that
$\{\overline{b}_{1},\cdots,\overline{b}_{n}\}$ generate the
simplicial fan $\Sigma$ (here $m\geq n$).
\end{defn}
Given an extended stacky fan $\mathbf{\Sigma}=(N,\Sigma,\beta)$, 
we have the following exact sequences:
\begin{equation}\Label{exact1}
0\longrightarrow DG(\beta)^{\star}\longrightarrow
\mathbb{Z}^{m}\stackrel{\beta}{\longrightarrow} N\longrightarrow
Coker(\beta)\longrightarrow 0,
\end{equation}
\begin{equation}\Label{exact2}
0\longrightarrow N^{\star}\longrightarrow
\mathbb{Z}^{m}\stackrel{\beta^{\vee}}{\longrightarrow}
DG(\beta)\longrightarrow Coker(\beta^{\vee})\longrightarrow 0,
\end{equation}
where $\beta^{\vee}$ is the Gale dual of $\beta$ (see \cite{BCS}). Applying $Hom_\mathbb{Z}(-,\mathbb{C}^{*})$ to (\ref{exact2}) yields
\begin{equation}\Label{exact3}
1\longrightarrow \mu\longrightarrow
G\stackrel{\alpha}{\longrightarrow}
(\mathbb{C}^{*})^{m}\longrightarrow
(\mathbb{C}^{*})^{d}\longrightarrow 1. \
\end{equation}
The toric 
DM stack $\mathcal{X}(\mathbf{\Sigma})$ is the quotient stack $[Z/G]$, where 
$Z:=(\mathbb{C}^{n}\setminus V(J_{\Sigma}))\times (\mathbb{C}^{*})^{m-n}$, $J_{\Sigma}$ is the 
irrelevant ideal of the fan $\Sigma$ and $G$ acts on $Z$ through the map $\alpha$ in (\ref{exact3}). The 
coarse moduli space of $\mathcal{X}(\mathbf{\Sigma})$ is the simplicial toric variety 
$X(\Sigma)$ corresponding to the simplicial fan $\Sigma$, see \cite{BCS} and \cite{Jiang}.

\begin{defn}\Label{semi-toric}
A toric DM stack $\mathcal{X}(\mathbf{\Sigma})$ is {\em semi-projective} if the coarse 
moduli space $X(\Sigma)$ is semi-projective.
\end{defn}

\begin{thm}\Label{semitor}
The following notions are equivalent:
\begin{enumerate}
\item A semi-projective toric DM stack $\mathcal{X}(\mathbf{\Sigma})$; 
\item A toric DM stack $\mathcal{X}(\mathbf{\Sigma})$ such that the simplicial fan $\Sigma$
is a regular triangulation of $\mathcal{B}=\{\overline{b}_{1},\cdots,\overline{b}_{n}\}$
which spans the lattice $\overline{N}$.
\end{enumerate}
\end{thm}

\begin{pf}
Since the toric DM stack is semi-projective if its coarse moduli space is semi-projective, 
the theorem follows from results in \cite{HS}.
\end{pf}

\subsection{The inertia stack}
Let $\mathbf{\Sigma}$ be an extended stacky fan and $\sigma\in\Sigma$ a cone. Define $link(\sigma):=\{\tau:\sigma+\tau\in \Sigma, \sigma\cap \tau=0\}$. Let $\{\widetilde{\rho}_{1},\ldots,\widetilde{\rho}_{l}\}$ be the rays in $link(\sigma)$. Consider the quotient extended stacky fan $\mathbf{\Sigma/\sigma}=(N(\sigma),\Sigma/\sigma,\beta(\sigma))$, with $\beta(\sigma): \mathbb{Z}^{l+m-n}\to N(\sigma)$ given by the images of $b_{1},\ldots,b_{l}$ and $b_{n+1},\ldots,b_{m}$
under $N\to N(\sigma)$. By the construction of toric Deligne-Mumford stacks, 
if $\sigma$ is contained in a top dimensional cone in $\Sigma$, we have $\mathcal{X}(\mathbf{\Sigma/\sigma}):=[Z(\sigma)/G(\sigma)]$, 
where $Z(\sigma)=(\mathbb{A}^{l}\setminus\mathbb{V}(J_{\Sigma/\sigma}))\times (\mathbb{C}^{*})^{m-n}$ and $G(\sigma)=Hom_{\mathbb{Z}}(DG(\beta(\sigma)),\mathbb{C}^{*})$.

\begin{lem}
If $\mathcal{X}(\mathbf{\Sigma})$ is semi-projective, so is $\mathcal{X}(\mathbf{\Sigma/\sigma})$.
\end{lem}

\begin{pf}
Semi-projectivity of the stack $\mathcal{X}(\mathbf{\Sigma})$ means the simplicial fan $\Sigma$ is a fan coming from a regular triangulation of $\mathcal{B}=\{\overline{b}_{1},\cdots,\overline{b}_{n}\}$ which spans the lattice $\overline{N}$. Let $pos(\mathcal{B})$ be the convex polyhedral cone generated by $\mathcal{B}$. Then from \cite{HS}, the triangulation is supported on $pos(\mathcal{B})$ and is dermined by a simple polyhedron whose normal fan is $\Sigma$. So $\sigma$ is contained in a top-dimensional 
cone $\tau$ in $\Sigma$. The image $\widetilde{\tau}$ of $\tau$ under quotient by $\sigma$ is a top-dimensional cone in the quotient fan $\Sigma/\sigma$. So the toric variety $X(\Sigma/\sigma)$ is semi-projective by Theorem \ref{semitor}, and the stack $\mathcal{X}(\mathbf{\Sigma/\sigma})$ is semi-projective by definition.
\end{pf}

Recall in \cite{BCS} that for each top-dimensional cone $\sigma$ in $\Sigma$, define $Box(\sigma)$ to be the set of elements $v\in N$ such that $\overline{v}=\sum_{\rho_{i}\subseteq \sigma}a_{i}\overline{b}_{i}$ for some $0\leq a_{i}<1$. Elements in  $Box(\sigma)$ are in one-to-one correspondence with elements in the finite group $N(\sigma)=N/N_{\sigma}$, where $N(\sigma)$ is a local group of the stack
$\mathcal{X}(\mathbf{\Sigma})$. In fact, we write $\overline{v}=\sum_{\rho_{i}\subseteq \sigma(\overline{v})}a_{i}\overline{b}_{i}$ for some $0<a_{i}<1$, where $\sigma(\overline{v})$ is the minimal cone containing $\overline{v}$. Denoted by $Box(\mathbf{\Sigma})$ the union  of $Box(\sigma)$ for all top-dimensional cones $\sigma$.  

\begin{prop}
The $r$-inertia stack is given by
\begin{equation}\Label{inertia}
\mathcal{I}_{r}\left(\mathcal{X}(\mathbf{\Sigma})\right)=\coprod_{(v_{1},\cdots,v_{r})\in Box(\mathbf{\Sigma})^{r}}
\mathcal{X}(\mathbf{\Sigma/\sigma}(\overline{v}_{1},\cdots,\overline{v}_{r})),
\end{equation}
where $\sigma(\overline{v}_{1},\cdots,\overline{v}_{r})$ is the minimal cone in 
$\Sigma$  containing $\overline{v}_{1},\cdots,\overline{v}_{r}$.
\end{prop} 
\begin{pf}
Since $G$ is an abelian group, we have 
$$\mathcal{I}_{r}\left(\mathcal{X}(\mathbf{\Sigma})\right)=[(\coprod_{(v_{1},\cdots,v_{r})\in (G)^{r}}
Z^{(v_{1},\cdots,v_{r})})\slash G],$$
where $Z^{(v_{1},\cdots,v_{r})}\subset Z$ is the subvariety fixed by $v_{1},\cdots,v_{r}$. Since 
$\sigma(\overline{v}_{1},\cdots,\overline{v}_{r})$ is contained in a top-dimensional cone  in 
$\Sigma$. We use the same method as in Lemma 4.6 and Proposition 4.7 of \cite{BCS} to prove that 
$[Z^{(v_{1},\cdots,v_{r})}\slash G]\cong \mathcal{X}(\mathbf{\Sigma/\sigma}(\overline{v}_{1},\cdots,\overline{v}_{r}))$.
\end{pf}

Note that in (\ref{inertia}) each component is semi-projective.

\subsection{The orbifold Chow ring}\label{ring}
In this section we compute the orbifold Chow ring of semi-projective toric DM stacks and prove Theorem \ref{main}.

\subsubsection{The module structure}

Let $\mathbf{\Sigma}=(N,\Sigma,\beta)$ be an extended  stacky fan such that the toric DM stack $\mathcal{X}(\mathbf{\Sigma})$ is semi-projective. Since the fan $\Sigma$ is convex, $|\Sigma|$ is an abelian subgroup of $N$. We put $N_{\Sigma}:=|\Sigma|\oplus N_{tor}$, where $N_{tor}$ is the torsion 
subgroup of $N$.  Define the deformed ring $\mathbb{Q}[N_{\Sigma}]:=\bigoplus_{c\in N_{\Sigma}}\mathbb{Q}y^{c}$ with the product structure given by (\ref{productA}).

Let $\{\rho_{1},\ldots,\rho_{n}\}$ be the rays of $\Sigma$, then each $\rho_{i}$ corresponds to a line
bundle $L_{i}$ over the toric Deligne-Mumford stack $\mathcal{X}(\mathbf{\Sigma})$ given by the trivial line bundle $\mathbb{C}\times Z$ over $Z$ with the $G$ action on $\mathbb{C}$ given by the $i$-th component $\alpha_{i}$ of $\alpha: G\to (\mathbb{C}^{*})^{m}$ in (\ref{exact3}). The first Chern classes of the line bundles $L_{i}$, which we identify with $y^{b_{i}}$, generate the cohomology ring of the simplicial toric variety $X(\Sigma)$.

Let $S_{\mathbf{\Sigma}}$ be the quotient ring $\frac{\mathbb{Q}[y^{b_{1}},\cdots,y^{b_{n}}]}{I_{\Sigma}}$, where $I_{\Sigma}$ is the square-free ideal of the fan $\Sigma$ generated by the monomials 
$$\{y^{b_{i_{1}}}\cdots y^{b_{i_{k}}}: \overline{b}_{i_{1}},\cdots, \overline{b}_{i_{k}} \text{ do not generate a cone in }\Sigma\}.$$
It is clear that $S_{\mathbf{\Sigma}}$ is a subring of the deformed ring $\mathbb{Q}[N_{\Sigma}]$.

\begin{lem}
Let $A^{*}(\mathcal{X}(\mathbf{\Sigma}))$ be the ordinary Chow ring of a semi-projective toric DM stack $\mathcal{X}(\mathbf{\Sigma})$. Then there is a ring isomorphism:
$$A^{*}(\mathcal{X}(\mathbf{\Sigma}))\cong \frac{S_{\mathbf{\Sigma}}}
{\{\sum_{i=1}^{n}e(b_{i})y^{b_{i}}: e\in N^{\star}\}}.$$
\end{lem}

\begin{pf}
The Lemma is easily proven from the fact that the Chow ring of a DM stack is isomorphic to the Chow ring of its coarse moduli space (\cite{V}) and Proposition 2.11 in \cite{HS}.
\end{pf}

Now we study the module structure on $A_{orb}^{*}\left(\mathcal{X}(\mathbf{\Sigma})\right)$. Because $\Sigma$ is a simplicial fan, we have:

\begin{lem}\Label{smalllemma}
For any $c\in N_{\Sigma}$, let $\sigma$ be the minimal cone in $\Sigma$ containing $\overline{c}$. Then there is a unique expression $c=v+\sum_{\rho_{i}\subset\sigma}m_{i}b_{i}$ where $m_{i}\in \mathbb{Z}_{\geq 0}$, and $v\in Box(\sigma)$.
\end{lem}

\begin{prop}\Label{vectorspace}
Let $\mathcal{X}(\mathbf{\Sigma})$ be a semi-projective toric  DM stack associated to an extended  stacky fan $\mathbf{\Sigma}$. We have an isomorphism of $A^{*}(\mathcal{X}(\mathbf{\Sigma}))$-modules:
$$\bigoplus_{v\in Box(\mathbf{\Sigma})}A^{*}\left(\mathcal{X}(\mathbf{\Sigma/\sigma}(\overline{v}))\right)[deg(y^{v})]\cong \frac{\mathbb{Q}[N_{\Sigma}]}{\{\sum_{i=1}^{n}e(b_{i})y^{b_{i}}: e\in N^{\star}\}}.$$
\end{prop}

\begin{pf}
From the definition of $\mathbb{Q}[N_{\Sigma}]$ and Lemma \ref{smalllemma}, we see  that $\mathbb{Q}[N_{\Sigma}]=\bigoplus_{v\in Box(\mathbf{\Sigma})}y^{v}\cdot S_{\mathbf{\Sigma}}$. The rest is similar to the proof of Proposition 4.7 in \cite{Jiang}, we leave it to the readers.
\end{pf}


\subsubsection{The Chen-Ruan product structure}

The orbifold cup product on a DM stack $\mathcal{X}$ is defined using genus zero, degree zero 3-pointed orbifold Gromov-Witten invariants on $\mathcal{X}$. The relevant moduli space is the disjoint union of all  3-twisted sectors (i.e. the double inertia stack). By (\ref{inertia}), the 3-twisted sectors of a semi-projective toric DM stack $\mathcal{X}(\mathbf{\Sigma})$ are 
\begin{equation}\Label{3-sector}
\coprod_{(v_{1},v_{2},v_{3})\in Box(\mathbf{\Sigma})^{3},
v_{1}v_{2}v_{3}=1}
~\mathcal{X}(\mathbf{\Sigma/\sigma}(\overline{v}_{1},\overline{v}_{2},\overline{v}_{3})).
\end{equation}

Let $ev_{i}:
\mathcal{X}(\mathbf{\Sigma/\sigma}(\overline{v}_{1},\overline{v}_{2},\overline{v}_{3}))\to
\mathcal{X}(\mathbf{\Sigma/\sigma}(\overline{v}_{i}))$ be the evaluation maps. The obstruction bundle (see \cite{CR2}) $Ob_{(v_{1},v_{2},v_{3})}$ over the 3-twisted sector
$\mathcal{X}(\mathbf{\Sigma/\sigma}(\overline{v}_{1},\overline{v}_{2},\overline{v}_{3}))$ are defined by 
\begin{equation}\Label{obstruction}
Ob_{(v_{1},v_{2},v_{3})}:=\left(e^{*}T\left(\mathcal{X}(\mathbf{\Sigma})\right)\otimes
H^{1}(C,\mathcal{O}_{C})\right)^{H},
\end{equation}
where $e: \mathcal{X}(\mathbf{\Sigma/\sigma}(\overline{v}_{1},\overline{v}_{2},\overline{v}_{3}))\to
\mathcal{X}(\mathbf{\Sigma})$ is the embedding,  $C\to \mathbb{P}^{1}$ is the $H$-covering branched over three marked points $\{0,1,\infty\}\subset \mathbb{P}^{1}$, and  $H$ is the group generated by $v_{1},v_{2},v_{3}$. 

A general result in \cite{CH} and \cite{JKK} about the obstruction bundle implies the following.

\begin{prop}\label{obstructionbdle}
Let 
$\mathcal{X}(\mathbf{\Sigma/\sigma}(\overline{v}_{1},\overline{v}_{2},\overline{v}_{3}))$ be a 3-twisted sector of the stack $\mathcal{X}(\mathbf{\Sigma})$. Suppose $v_{1}+v_{2}+v_{3}=\sum_{\rho_{i}\subset
\sigma(\overline{v}_{1},\overline{v}_{2},\overline{v}_{3})}a_{i}b_{i}$, $a_{i}=1$ or $2$. Then the Euler class of the obstruction bundle $Ob_{(v_{1},v_{2},v_{3})}$ on
$\mathcal{X}(\mathbf{\Sigma/\sigma}(\overline{v}_{1},\overline{v}_{2},\overline{v}_{3}))$
is
$$\prod_{a_{i}=2}c_{1}(L_{i})|_{\mathcal{X}(\mathbf{\Sigma/\sigma}(\overline{v}_{1},\overline{v}_{2},\overline{v}_{3}))},$$
where $L_{i}$ is the line bundle over $\mathcal{X}(\mathbf{\Sigma})$ corresponding to the ray $\rho_{i}$.
\end{prop}

Let $v\in Box(\mathbf{\Sigma})$, say $v\in N(\sigma)$ for some top-dimensional cone $\sigma$. Let $\check{v}\in Box(\mathbf{\Sigma})$ be the inverse of $v$ as an element in the group $N(\sigma)$. Equivalently, if $v=\sum_{\rho_{i}\subseteq \sigma(\overline{v})}\alpha_{i}b_{i}$ for $0<\alpha_{i}<1$, then $\check{v}=\sum_{\rho_{i}\subseteq \sigma(\overline{v})}(1-\alpha_{i})b_{i}$. Then for $\alpha_{1},\alpha_{2}\in A^{*}_{orb}(\mathcal{X}(\mathbf{\Sigma}))$, the orbifold cup product is defined by
\begin{equation}\Label{cupproduct} 
\alpha_{1}\cup_{orb}\alpha_{2}=\widehat{ev}_{3*}(ev_{1}^{*}\alpha_{1}\cup ev_{2}^{*}\alpha_{2}\cup e(Ob_{(v_{1},v_{2},v_{3})})),
\end{equation}
where $\widehat{ev}_{3}=I\circ ev_{3}$, and $I: \mathcal{I}\mathcal{X}(\mathbf{\Sigma}) \rightarrow \mathcal{I}\mathcal{X}(\mathbf{\Sigma})$ is the natural map given by  $(x,g)\mapsto (x,g^{-1})$.

\subsubsection*{Proof of Theorem 1.1}
By Proposition \ref{vectorspace}, it remains to consider the cup product. In this case, for any 
$v_{1},v_{2}\in Box(\mathbf{\Sigma})$, we also have
$$v_{1}+v_{2}=\check{v}_{3}+\sum_{a_{i}=2}b_{i}+\sum_{i\in J}b_{i},$$
where $J$ represents the set of $j$
such that $\rho_{j}$ belongs to
$\sigma(\overline{v}_{1},\overline{v}_{2})$, but not belong to
$\sigma(\overline{v}_{3})$. Then the proof is the same as the proof in \cite{BCS}. 
We omit the details.

\section{Lawrence Toric DM stacks}\label{Lawrence}

In this section we study a special type of semi-projective toric DM stacks called the Lawrence toric DM stacks. Their orbifold Chow rings are shown to be isomorphic to the orbifold Chow rings of their associated hypertoric DM stacks studied in \cite{JT}.

\subsection{Stacky hyperplane arrangements}
Let $N$, $\{b_{1},\cdots,b_{m}\}\in N$, $\beta:\mathbb{Z}^{m}\to N$, and $\{\overline{b}_{1},\cdots,\overline{b}_{m}\}\subset \overline{N}$ be as in Definition \ref{stackyfan}.  We assume that $\{b_{1},\cdots,b_{m}\}\in N$ are nontorsion integral vectors. We still have the exact sequences (\ref{exact1}) and (\ref{exact2}). The Gale dual map $\beta^{\vee}$ of $\beta$ is given by a collection of integral vectors $\beta^{\vee}=(a_1,\cdots,a_m)$. Choose a generic element $\theta\in DG(\beta)$ and let $\psi:=(r_{1},\cdots,r_{m})$ be a lifting of $\theta$ in $\mathbb{Z}^{m}$ such that $\theta=-\beta^{\vee}\psi$. Note that $\theta$ is generic if and only if it is not in any hyperplane of the configuration determined by $\beta^{\vee}$ in $DG(\beta)_{\mathbb{R}}$. Associated to $\theta$ there is a hyperplane arrangement $\mathcal{H}=\{H_{1},\cdots,H_{m}\}$ defined as follows:  let $H_{i}$ be the hyperplane
\begin{equation}\Label{arrangement}
H_{i}:=\{v\in M_{\mathbb{R}}|<\overline{b}_{i},v>+r_{i}=0\}\subset M_{\mathbb{R}}.
\end{equation}
This determines hyperplane arrangement in $M_{\mathbb{R}}$, up to translation. 
It is well-known that hyperplane arrangements determine the topology of hypertoric varieties (\cite{BD}). 
We call $\mathcal{A}:=(N,\beta,\theta)$ a {\em stacky hyperplane arrangement}.

The toric variety $X(\Sigma)$ is defined by the weighted polytope $\mathbf{\Gamma}:=\bigcap_{i=1}^{m}F_{i}$, where $F_{i}=\{v\in M_{\mathbb{R}}|<b_{i},v>+r_{i}\geq 0\}$. Suppose that $\mathbf{\Gamma}$ is bounded, the  fan $\Sigma$ is the normal fan of  $\mathbf{\Gamma}$ in $M_{\mathbb{R}}=\mathbb{R}^{d}$ with one dimensional rays generated by $\overline{b}_{1},\cdots,\overline{b}_{n}$. By reordering, we may assume that $H_{1},\cdots,H_{n}$ are the hyperplanes that bound the polytope $\mathbf{\Gamma}$, and $H_{n+1},\cdots,H_{m}$ are the other hyperplanes. Then we have an extended stacky fan
$\mathbf{\Sigma}=(N,\Sigma,\beta)$ as in Definition \ref{stackyfan}, with $\Sigma$ the normal fan of $\mathbf{\Gamma}$, $\beta:\mathbb{Z}^{m}\to N$ given by $\{b_{1},\cdots,b_{n},b_{n+1},\cdots,b_{m}\}\subset N$, and $\{b_{n+1},\cdots,b_{m}\}$ the extra data. We define the hypertoric DM stack $\mathcal{M}(\mathcal{A})$ using this $\mathcal{A}$, see \cite{JT} for more details.

\subsection{Lawrence toric DM stacks}
Applying Gale dual to the map 
\begin{equation}\label{betaL}
\mathbb{Z}^{m}\oplus \mathbb{Z}^{m}\to DG(\beta),
\end{equation}
given by $(\beta^{\vee},-\beta^{\vee})$, we obtain
$$\beta_{L}: \mathbb{Z}^{m}\oplus \mathbb{Z}^{m}\longrightarrow N_{L},$$
which  is given by integral  vectors $\{b_{L,1},\cdots,b_{L,m},b'_{L,1},\cdots,b'_{L,m}\}$ in $N_{L}$. The natural images  $\{\overline{b}_{L,1},\cdots,\overline{b}_{L,m},\overline{b}'_{L,1},\cdots,\overline{b}'_{L,m}\}\subset \overline{N}_{L}$ 
are called the Lawrence lifting of $\{\overline{b}_{1},\cdots,\overline{b}_{m}\}\subset \overline{N}$.  

Associated to the generic element $\theta$, let $\overline{\theta}$ be the natural image under the map $DG(\beta)\rightarrow \overline{DG(\beta)}$. Then the map $\overline{\beta}^{\vee}: \mathbb{Z}^{m}\rightarrow \overline{DG(\beta)}$ is given by $\overline{\beta}^{\vee}=(\overline{a}_1,\cdots,\overline{a}_m)$. For any column basis of the form $C=\{\overline{a}_{i_{1}},\cdots,\overline{a}_{i_{m-d}}\}$, there exist unique $\lambda_{1},\cdots,\lambda_{m-d}$ such that
$$a_{i_{1}}\lambda_{1}+\cdots+a_{i_{m-d}}\lambda_{m-d}=\overline{\theta}.$$
Let $\mathbb{C}[z_{1},\cdots,z_{m},w_{1},\cdots,w_{m}]$ be the coordinate ring of $\mathbb{C}^{2m}$. Let
$\sigma(C,\theta)=\{\overline{b}_{i_{j}}~|\lambda_{j}>0\}\sqcup\{\overline{b}'_{i_{j}}|~\lambda_{j}<0\},$
and 
$C(\theta)=\{z_{i_{j}}~|\lambda_{j}>0\}\sqcup\{w_{i_{j}}|~\lambda_{j}<0\}$.
We set
\begin{equation}\Label{irrelevant}
\mathbf{\mathcal{I}}_{\theta}:=<\prod
C(\theta)|~C~\text{is a column basis of}~\overline{\beta}^{\vee}>,
\end{equation}
and 
\begin{equation}\Label{fan}
\Sigma_{\theta}:=\{\overline{\sigma}(C,\theta):~C~\text{is a column basis of}~\overline{\beta}^{\vee}\},
\end{equation}
where $\overline{\sigma}(C,\theta)=
\{\overline{b}_{L,1},\cdots,\overline{b}_{L,m},\overline{b}'_{L,1},\cdots,\overline{b}'_{L,m}\}\setminus\sigma(C,\theta)$ 
is the complement of $\sigma(C,\theta)$ and corresponds to the maximal cones in $\Sigma_{\theta}$.
According to \cite{HS}, $\Sigma_{\theta}$ is the fan  of
Lawrence toric variety $X(\Sigma_{\theta})$ corresponding to
$\theta$ in the lattice $\overline{N}_{L}$. 
The ideal
$\mathcal{I}_{\theta}$ is the irrelevant ideal of the fan $\Sigma_{\theta}$. 
Then we have the Lawrence stacky fan $\mathbf{\Sigma_{\theta}}=(N_{L},\Sigma_{\theta},\beta_{L})$ introduced in \cite{JT}.

Applying $Hom_\mathbb{Z}(-,\mathbb{C}^{*})$ functor to (\ref{betaL}), we get 
\begin{equation}\Label{Lawrencemap}
\alpha_{h}: G\rightarrow (\mathbb{C}^{*})^{2m}.
\end{equation}
So $G$ acts on $\mathbb{C}^{2m}$ through $\alpha_{h}$. From Section 2,  
$\mathcal{X}(\mathbf{\Sigma_{\theta}})=[(\mathbb{C}^{2m}\setminus V(\mathcal{I}_{\theta}))\slash G]$ whose coarse moduli space is the Lawrence toric variety $X(\Sigma_{\theta})=(\mathbb{C}^{2m}\setminus V(\mathcal{I}_{\theta}))\slash G$. Let $Y\subset \mathbb{C}^{2m}\setminus V(\mathcal{I}_{\theta})$
be the subvariety defined by the ideal:
\begin{equation}\Label{ideal1}
I_{\beta^{\vee}}:=<\sum_{i=1}^{m}(\beta^{\vee})^{\star}(x)_{i}a_{ij}z_{i}w_{i}|\forall x\in DG(\beta)^{\star}>,
\end{equation}
where $(\beta^{\vee})^{\star}: DG(\beta)^{\star}\rightarrow \mathbb{Z}^{m}$ is the dual map of $\beta^{\vee}$ and $(\beta^{\vee})^{\star}(x)_{i}$ is the $i$-th component
of the vector $(\beta^{\vee})^{\star}(x)$.
From \cite{JT}, the hypertoric DM stack $\mathcal{M}(\mathcal{A})=[Y/G]$ whose coarse moduli space is the hypertoric variety $Y(\beta^{\vee},\theta)=Y\slash G$.

\begin{defn}(\cite{JT})
The Lawrence toric DM stack is  the toric DM stack $\mathcal{X}(\mathbf{\Sigma_{\theta}})$ corresponding to the Lawrence stacky fan $\mathbf{\Sigma_{\theta}}$.
\end{defn}

By \cite{HS},  $X(\Sigma_{\theta})$ is semi-projective. So the Lawrence toric DM stack $\mathcal{X}(\mathbf{\Sigma_{\theta}})$ is semi-projective by definition.

\subsection{Comparison of inertia stacks}
Next we compare the orbifold Chow ring of the hypertoric DM stack and the orbifold Chow ring of the Lawrence toric DM stack. First we compare the inertia stacks. From the map 
$\beta: \mathbb{Z}^{m}\rightarrow N$ which is given by vectors $\{b_{1},\cdots,b_m\}$.  Let $Cone(\beta)$ be a partially ordered finite set of cones generated by $\overline{b}_{1},\cdots,\overline{b}_{m}$. The partial order is defined by: $\sigma\prec\tau$ if $\sigma$ is a face of $\tau$, and we have the minimum element $\hat{0}$ which is the cone consisting of the origin. Let $Cone(\overline{N})$ be the set of all convex polyhedral cones in the lattice $\overline{N}$. Then we have a map
$$C: Cone(\beta)\longrightarrow Cone(\overline{N}),$$ such that for any $\sigma\in Cone(\beta)$, $C(\sigma)$ is the cone in $\overline{N}$. Then $\Delta_{\mathbf{\beta}}:=(C,Cone(\beta))$ is a simplicial {\em multi-fan} in the sense of \cite{HM}. 

For the multi-fan $\Delta_{\mathbf{\beta}}$, let  $Box(\Delta_{\mathbf{\beta}})$ be the set of pairs $(v,\sigma)$, where $\sigma$ is a cone in $\Delta_{\mathbf{\beta}}$, $v\in N$ such that $\overline{v}=\sum_{\rho_{i}\subseteq \sigma}\alpha_{i}b_{i}$ for $0<\alpha_{i}<1$. (Note that $\sigma$ is the minimal cone in $\Delta_{\mathbf{\beta}}$ satisfying the above condition.)  From \cite{JT}, an element $(v,\sigma)\in Box(\Delta_{\mathbf{\beta}})$ gives a component of the inertia stack $\mathcal{I}(\mathcal{M}(\mathcal{A}))$.
Also consider the set $Box(\mathbf{\Sigma_{\theta}})$ associated to the stacky fan $\mathbf{\Sigma_{\theta}}$, see Section 2.2 for its definition. An element $v\in Box(\mathbf{\Sigma_{\theta}})$ gives a component of the inertia stack $\mathcal{I}(\mathcal{X}(\mathbf{\Sigma_{\theta}}))$. 

By the Lawrence lifting property, a vector $\overline{b}_{i}$ in $\overline{N}$ lifts to two vectors
$\overline{b}_{L,i},\overline{b}'_{L,i}$ in $\overline{N}_{L}$. Let $\{\overline{b}_{L,i_{1}},\cdots,\overline{b}_{L,i_{k}}, \overline{b}'_{L,i_{1}},\cdots,\overline{b}'_{L,i_{k}}\}$ be the Lawrence lifting of $\{\overline{b}_{i_{1}},\cdots,\overline{b}_{i_{k}}\}$.

\begin{lem}\Label{conemulti}
$\{\overline{b}_{i_{1}},\cdots,\overline{b}_{i_{k}}\}$ generate a cone  $\sigma$ in $\Delta_{\mathbf{\beta}}$ if and only if $\{\overline{b}_{L,i_{1}},\cdots,\overline{b}_{L,i_{k}}, \overline{b}'_{L,i_{1}},\cdots,\overline{b}'_{L,i_{k}}\}$  generate a cone $\sigma_{\theta}$ in $\Sigma_{\theta}$.
\end{lem}

\begin{pf}
Suppose $\sigma$ is a cone in $\Delta_{\mathbf{\beta}}$ generated by $\{\overline{b}_{i_{1}},\cdots,\overline{b}_{i_{k}}\}$, it is contained in a top-dimensional cone $\tau$. Assume that $\tau$ is generated by $\{\overline{b}_{i_{1}},\cdots,\overline{b}_{i_{k}}, \overline{b}_{i_{k+1}},\cdots,\overline{b}_{i_{d}}\}$. Let $C$ be the complement $\{\overline{b}_{1},\cdots,\overline{b}_{m}\}\setminus \tau$. Then $C$ corresponds to a column basis of  $\overline{\beta}^{\vee}$ in the map $\overline{\beta}^{\vee}: \mathbb{Z}^{m}\rightarrow \overline{DG(\beta)}$. By the definition of $\Sigma_{\theta}$ in (\ref{fan}), $C$ corresponds to a maximal cone $\tau_{\theta}$ in $\Sigma_{\theta}$ which contains
the rays generated by $\{\overline{b}_{L,i_{1}},\cdots,\overline{b}_{L,i_{k}},
\overline{b}'_{L,i_{1}},\cdots,\overline{b}'_{L,i_{k}}\}$. Thus these rays generate a cone $\sigma_{\theta}$
in $\Sigma_{\theta}$.  

Conversely, suppose $\sigma_{\theta}$ is a cone in $\Sigma_{\theta}$ generated by $\{\overline{b}_{L,i_{1}},\cdots,\overline{b}_{L,i_{k}}, \overline{b}'_{L,i_{1}},\cdots,\overline{b}'_{L,i_{k}}\}$. Using the similar method above we prove that $\{\overline{b}_{i_{1}},\cdots,\overline{b}_{i_{k}}\}$ must be contained in a top-dimensional cone of $\Delta_{\mathbf{\beta}}$. So $\{\overline{b}_{i_{1}},\cdots,\overline{b}_{i_{k}}\}$ generate a cone $\sigma$ in $\Delta_{\mathbf{\beta}}$.
\end{pf}

\begin{lem}\Label{box}
There is an one-to-one correspondence between the elements in $Box(\mathbf{\Sigma_{\theta}})$ and the elements in $Box(\Delta_{\mathbf{\beta}})$. Moreover, their degree shifting numbers coincide.
\end{lem}

\begin{pf}
First the torsion elements in $Box(\mathbf{\Sigma_{\theta}})$ and $Box(\Delta_{\mathbf{\beta}})$ are both 
isomorphic to $\mu=ker(\alpha)=ker(\alpha_{h})$ in (\ref{exact3}) and (\ref{Lawrencemap}).
Let $(v,\sigma)\in Box(\Delta_{\mathbf{\beta}})$ with $\overline{v}=\sum_{\rho_i\subseteq \sigma}\alpha_{i}\overline{b}_{i}$.
Then $v$ may be identified with an element (which we ambiguously denote by) $v\in G:=Hom_\mathbb{Z}(DG(\beta),\mathbb{C}^*)$. Certainly $v$ fixes a point in $\mathbb{C}^{m}$. Consider the map $\alpha$ in (\ref{exact3}), put $\alpha(v)=(\alpha^{1}(v),\cdots,\alpha^{m}(v))$. Then $\alpha^{i}(v)\neq 1$ if $\rho_i\subseteq\sigma$, and $\alpha^{i}(v)= 1$ otherwise. By Lemma \ref{conemulti}, let $\{\overline{b}_{L,i},\overline{b}'_{L,i}:i=1,\cdots,|\sigma|\}$ be the Lawrence lifting of $\{\overline{b}_{i}\}_{\rho_i\subseteq\sigma}$. Since the action of $v$ on $\mathbb{C}^{2m}$ is given by $(v,v^{-1})$, $v$ fixes a point in $\mathbb{C}^{2m}$ and yields an element $v_\theta$ in $Box(\mathbf{\Sigma_\theta})$. 
From the map (\ref{Lawrencemap}), let
\begin{equation}\Label{vtheta}
\alpha_{h}(v_{\theta})=(\alpha^{1}_{h}(v_{\theta}),\cdots,\alpha^{m}_{h}(v_{\theta}),
\alpha^{m+1}_{h}(v_{\theta}),\cdots,\alpha^{2m}_{h}(v_{\theta})).
\end{equation}
Then $\alpha_{h}^{i}(v_{\theta})\neq 1$ and $\alpha_{h}^{i+m}(v_{\theta})\neq 1$ if $\rho_i\subseteq\sigma$; $\alpha_{h}^{i}(v_{\theta})= \alpha_{h}^{i+m}(v_{\theta})= 1$ otherwise.  So $\sigma_{\theta}(\overline{v}_{\theta})=\{\overline{b}_{L,i},\overline{b}'_{L,i}:i=1,\cdots,|\sigma|\}$ is the minimal cone in $\Sigma_{\theta}$ containing $\overline{v}_{\theta}$. Furthermore, $\overline{v}_{\theta}=\sum_{\rho_i\subseteq \sigma}\alpha_{i}\overline{b}_{L,i}+ \sum_{\rho_i\subseteq \sigma}(1-\alpha_{i})\overline{b}'_{L,i}$.

Conservely, given an element $v_{\theta}\in Box(\mathbf{\Sigma_{\theta}})$, let $\sigma_{\theta}(\overline{v}_{\theta})$ be the minimal cone in $\Sigma_{\theta}$ containing $\overline{v}_{\theta}$. 
Then from the action of $G$ on $\mathbb{C}^{2m}$ and (\ref{vtheta}), we have $\alpha^{i}_{h}(v_{\theta})=(\alpha^{i+m}_{h}(v_{\theta}))^{-1}$. If $\alpha^{i}_{h}(v_{\theta})\neq 1$, then $\alpha^{i+m}_{h}(v_{\theta})\neq 1$, which means that  $\overline{b}_{L,i}, \overline{b}_{L,i+m}\in \sigma_{\theta}(\overline{v}_{\theta})$. The cone $\sigma_{\theta}(\overline{v}_{\theta})$ is the one in $\Sigma_{\theta}$ containing $\overline{b}_{L,i}, \overline{b}_{L,i+m}$'s satisfying this condition. 
Then $\overline{v}_{\theta}=\sum_{i}(\alpha_{i}\overline{b}_{L,i}+(1-\alpha_{i})\overline{b}^{'}_{L,i})$. 
By Lemma \ref{conemulti},
$\sigma_{\theta}(\overline{v}_{\theta})$ is the Lawrence lifting of a cone $\sigma$ generated by the $\{\overline{b}_{i}\}$'s in $\Delta_{\mathbf{\beta}}$. Let $v=\sum_{\rho_i\subseteq\sigma}\alpha_{i}b_{i}$. So it also determines an element $(v,\sigma)\in Box(\Delta_{\mathbf{\beta}})$.
\end{pf}

For $(v_{1},\sigma_{1}),(v_{2},\sigma_{2}),(v_{3},\sigma_{3})\in Box(\Delta_{\mathbf{\beta}})$, let $\sigma(\overline{v}_{1}, \overline{v}_{2},\overline{v}_{3})$ be the miniaml cone containing $\overline{v}_{1}, \overline{v}_{2},\overline{v}_{3}$ in $\Delta_{\mathbf{\beta}}$ such that $\overline{v}_{1}+\overline{v}_{2}+\overline{v}_{3}=\sum_{\rho_i\subseteq \sigma(\overline{v}_{1}, \overline{v}_{2},\overline{v}_{3})}a_{i}\overline{b}_{i}$ and $a_{i}=1,2$. Let $v_{\theta,1}, v_{\theta,2},v_{\theta,3}$ be the corresponding elements in  $Box(\mathbf{\Sigma_{\theta}})$ and $\sigma(\overline{v}_{\theta,1}, \overline{v}_{\theta,2},\overline{v}_{\theta,3})$ the minimal cone containing $\overline{v}_{\theta,1}, \overline{v}_{\theta,2},\overline{v}_{\theta,3}$ in $\mathbf{\Sigma_{\theta}}$. Then by Lemmas \ref{conemulti} and \ref{box}, 
$\sigma(\overline{v}_{\theta,1}, \overline{v}_{\theta,2},\overline{v}_{\theta,3})$ is the Lawrence lifting of  $\sigma(\overline{v}_{1}, \overline{v}_{2},\overline{v}_{3})$. Suppose that $\sigma$
is generated by $\{\overline{b}_{i_{1}},\cdots,\overline{b}_{i_{s}}\}$, then  $\sigma(\overline{v}_{\theta,1}, \overline{v}_{\theta,2},\overline{v}_{\theta,3})$ is generated by $\{\overline{b}_{L,i_{1}},\cdots,\overline{b}_{L,i_{s}},\overline{b}^{'}_{L,i_{1}},\cdots,\overline{b}^{'}_{L,i_{s}}\}$, the Lawrence lifting
of $\{\overline{b}_{i_{1}},\cdots,\overline{b}_{i_{s}}\}$. Let $\{\overline{b}_{j_{1}},\cdots,\overline{b}_{j_{m-l-s}}\}$ be the rays
not in $\sigma\cup link(\sigma)$, we have the Lawrence lifting 
$\{\overline{b}_{L,j_{1}},\cdots,\overline{b}_{L,j_{m-l-s}},\overline{b}^{'}_{L,j_{1}},\cdots,\overline{b}^{'}_{L,j_{m-l-s}}\}$. Then 
from the definition of Lawrence fan $\Sigma_{\theta}$ in (\ref{fan}), we have the following lemma:
\begin{lem}\Label{keycone}
There exist $m-l-s$  vectors in $\{\overline{b}_{L,j_{1}},\cdots,\overline{b}_{L,j_{m-l-s}},\overline{b}^{'}_{L,j_{1}},\cdots,\overline{b}^{'}_{L,j_{m-l-s}}\}$
such that the rays they generate plus the rays in $\sigma(\overline{v}_{\theta,1}, \overline{v}_{\theta,2},\overline{v}_{\theta,3})$ 
generate a cone $\sigma_{\theta}$ in $\Sigma_{\theta}$. $\square$
\end{lem}

\begin{prop}\Label{3-twisted-sector}
The stack $\mathcal{X}(\mathbf{\Sigma_{\theta}}/\sigma_{\theta})$ is also a Lawrence toric DM stack.
\end{prop}

\begin{pf}
For simplicity, put $\sigma:=\sigma(\overline{v}_{1},
\overline{v}_{2},\overline{v}_{3})$. Suppose there are
$l$ rays in the $link(\sigma)$. Then by Lemma \ref{conemulti} there
are $2l$ rays in $link(\sigma_{\theta})$, the Lawrence lifting of $link(\sigma)$.
Let $s:=|\sigma|$, then $2s+m-l-s=|\sigma_{\theta}|$. Applying Gale dual to the diagrams
\[
\begin{CD}
0 @ >>>\mathbb{Z}^{s}@ >>> \mathbb{Z}^{l+s}@ >>> \mathbb{Z}^{l} @
>>> 0\\
&& @VV{\beta_{\sigma}}V@VV{\widetilde{\beta}}V@VV{\beta(\sigma)}V \\
0@ >>>N_{\sigma} @ >{}>>N@ >>> N(\sigma) @>>> 0,
\end{CD}
\]
and 
\[
\begin{CD}
0 @ >>>\mathbb{Z}^{l+s}@ >>> \mathbb{Z}^{m}@ >>> \mathbb{Z}^{m-l-s} @
>>> 0\\
&& @VV{\widetilde{\beta}}V@VV{\beta}V@VV{}V \\
0@ >>>N @ >{\cong}>>N@ >>> 0 @>>> 0
\end{CD}
\]
yields
\begin{equation}\Label{3-sector2}
\begin{CD}
0 @ >>>\mathbb{Z}^{l}@ >>> \mathbb{Z}^{l+s}@ >>> \mathbb{Z}^{s} @
>>> 0\\
&& @VV{\beta(\sigma)^{\vee}}V@VV{\widetilde{\beta}^{\vee}}V@VV{\beta_{\sigma}^{\vee}}V \\
0@ >>>DG(\beta(\sigma)) @ >{\varphi_{1}}>>DG(\widetilde{\beta})@ >>> DG(\beta_{\sigma})
@>>> 0,
\end{CD}
\end{equation}
and 
\begin{equation}\Label{3-sector22}
\begin{CD}
0 @ >>>\mathbb{Z}^{m-l-s}@ >>> \mathbb{Z}^{m}@ >>> \mathbb{Z}^{l+s} @
>>> 0\\
&& @VV{\cong}V@VV{\beta^{\vee}}V@VV{\widetilde{\beta}^{\vee}}V \\
0@ >>>\mathbb{Z}^{m-l-s} @ >{}>>DG(\beta)@ >{\varphi_{2}}>> DG(\widetilde{\beta})
@>>> 0.
\end{CD}
\end{equation}
Since $\mathbb{Z}^{s}\cong N_{\sigma}$, we have that $DG(\beta_{\sigma})=0$. 
We add two exact sequences
$$0\longrightarrow \mathbb{Z}^{l}\longrightarrow\mathbb{Z}^{m}\longrightarrow\mathbb{Z}^{m-l}\longrightarrow 0,$$
and 
$$0\longrightarrow 0\longrightarrow\mathbb{Z}^{m}\longrightarrow\mathbb{Z}^{m}\longrightarrow 0,$$ 
on the rows of the diagrams (\ref{3-sector2}),(\ref{3-sector22}) and make suitable maps to  the
Gale duals we get
\begin{equation}\Label{3-sectors}
\begin{CD}
0 @ >>>\mathbb{Z}^{2l}@ >>> \mathbb{Z}^{l+s+m}@ >>>
\mathbb{Z}^{s+m-l} @
>>> 0\\
&& @VV{(\beta(\sigma)^{\vee},-\beta(\sigma)^{\vee})}V@VV{(\widetilde{\beta}^{\vee},-\beta^{\vee})}V@VV{0}V \\
0@ >>>DG(\beta(\sigma)) @ >{\cong}>>DG(\widetilde{\beta})@ >>> 0
@>>> 0,
\end{CD}
\end{equation}
and 
\begin{equation}\Label{3-sectors2}
\begin{CD}
0 @ >>>\mathbb{Z}^{m-l-s}@ >>> \mathbb{Z}^{2m}@ >>>
\mathbb{Z}^{l+s+m} @
>>> 0\\
&& @VV{\cong}V@VV{(\beta^{\vee},-\beta^{\vee})}V@VV{(\widetilde{\beta}^{\vee},-\beta^{\vee})}V \\
0@ >>>\mathbb{Z}^{m-l-s} @ >{}>>DG(\beta)@ >>> DG(\widetilde{\beta})
@>>> 0.
\end{CD}
\end{equation}
Applying Gale dual to (\ref{3-sectors}), (\ref{3-sectors2}) we get
\[
\begin{CD}
0 @ >>>\mathbb{Z}^{s+m-l}@ >>> \mathbb{Z}^{l+s+m}@ >>> \mathbb{Z}^{2l} @ >>> 0\\
&& @VV{\cong}V@VV{\widetilde{\beta}_{L}}V@VV{\beta_{L}(\sigma_{\theta})}V \\
0@ >>>\mathbb{Z}^{s+m-l} @ >{}>>\widetilde{N}_{L}@ >>> N_{L}(\sigma_{\theta})
@>>> 0,
\end{CD}
\]
 and
\[
\begin{CD}
0 @ >>>\mathbb{Z}^{l+s+m}@ >>> \mathbb{Z}^{2m}@ >>> \mathbb{Z}^{m-l-s} @ >>> 0\\
&& @VV{\widetilde{\beta}_{L}}V@VV{\beta_{L}}V@VV{0}V \\
0@ >>>\widetilde{N}_{L} @ >{\cong}>>N_{L}@ >>> 0
@>>> 0.
\end{CD}
\] 
For the generic element 
$\theta$, from them map $\varphi_{2}$ in (\ref{3-sector22}), 
$\theta$ induces $\widetilde{\theta}\in DG(\widetilde{\beta})$, and from the isomorphism
$\varphi_{1}$ in (\ref{3-sector2}),  $\widetilde{\theta}=\theta(\sigma)\in DG(\beta(\sigma))$.
So we a quotient stacky hyperplane arrangement $\mathcal{A}(\sigma)=(N(\sigma),\beta(\sigma),\theta(\sigma))$.
From the above diagrams  we see that the quotient fan
$\Sigma_{\theta}/\sigma_{\theta}$ in
$\overline{N}_{L}(\sigma_{\theta})$ also comes from a Lawrence
construction of the map $\beta(\sigma)^{\vee}:
\mathbb{Z}^{l}\rightarrow DG(\beta(\sigma))$.  Let
$X(\sigma)=\mathbb{C}^{2l}\setminus
V(\mathcal{I}_{\theta(\sigma)})$, where
$\mathcal{I}_{\theta(\sigma)}$ is the irrelevant ideal  of the
quotient fan $\Sigma_{\theta}\slash \sigma_{\theta}$. Let
$G(\sigma)=Hom_\mathbb{Z}(DG(\beta(\sigma)),\mathbb{C}^{*})$. The
stack
$\mathcal{X}(\mathbf{\Sigma_{\theta}}/\sigma_{\theta})=[X(\sigma)/G(\sigma)]$
is a Lawrence toric Deligne-Mumford stack.
\end{pf}

\begin{cor}\label{sectors}
$\mathcal{M}(\mathcal{A}(\sigma(\overline{v}_{1},
\overline{v}_{2},\overline{v}_{3})))$ is the hypertoric DM stack associated to the quotient Lawrence toric DM stack
$\mathcal{X}(\mathbf{\Sigma_{\theta}}/\sigma_{\theta})$.
\end{cor}

\begin{pf}
$\mathcal{M}(\mathcal{A}(\sigma(\overline{v}_{1},
\overline{v}_{2},\overline{v}_{3})))$ is constructed in \cite{JT} as a quotient stack 
$[Y(\sigma)/G(\sigma)]$, where $Y(\sigma)\subset X(\sigma)$ is defined by $I_{\beta(\sigma)^{\vee}}$, which is 
the ideal in (\ref{ideal1}) corresponding to the map $\beta(\sigma)^{\vee}$ in (\ref{3-sector2}). So the 
stack $\mathcal{M}(\mathcal{A}(\sigma(\overline{v}_{1},
\overline{v}_{2},\overline{v}_{3})))$ is the associated hypertoric DM stack in 
the Lawrence toric DM stack $\mathcal{X}(\mathbf{\Sigma_{\theta}}/\sigma_{\theta})$.
\end{pf}

\begin{rmk}
For any $v_{\theta}\in Box(\mathbf{\Sigma_{\theta}})$, let $v^{-1}_{\theta}$ be its inverse. We have
the quotient Lawrence toric stack $\mathcal{X}(\mathbf{\Sigma_{\theta}}/\sigma_{\theta})$. 
Let $(v,\sigma)$ be the corresponding element in $Box(\Delta_{\mathbf{\beta}})$, then  
$$\mathcal{M}(\mathcal{A}(\sigma(\overline{v},
\overline{v}^{-1},1)))\cong \mathcal{M}(\mathcal{A}(\sigma)).$$
By Proposition \ref{3-twisted-sector} and Corollary \ref{sectors}, the twisted sector
$\mathcal{M}(\mathcal{A}(\sigma))$ is the associated hypertoric DM stack 
of the Lawrence toric DM stack $\mathcal{X}(\mathbf{\Sigma_{\theta}}/\sigma_{\theta})$.
\end{rmk}

\begin{rmk}
From Lemma \ref{keycone}, the cone $\sigma_{\theta}$ is not the minimal cone 
$\sigma(\overline{v}_{\theta,1}, \overline{v}_{\theta,2},\overline{v}_{\theta,3})$ containing 
$\overline{v}_{\theta,1}, \overline{v}_{\theta,2},\overline{v}_{\theta,3}$ in $\Sigma_{\theta}$.
So $\mathcal{X}(\Sigma_{\theta}/\sigma(\overline{v}_{\theta,1}, \overline{v}_{\theta,2},\overline{v}_{\theta,3}))$
is not a Lawrence toric DM stack. But from the construction of Lawrence toric DM stack,
the quotient stack $\mathcal{X}(\Sigma_{\theta}/\sigma(\overline{v}_{\theta,1}, \overline{v}_{\theta,2},\overline{v}_{\theta,3}))$
is homotopy equivalent to the quotient stack $\mathcal{X}(\Sigma_{\theta}/\sigma_{\theta})$.
Since we do not need this to compare the orbifold Chow ring, we omit the details.
\end{rmk}

\subsection{Comparison of orbifold Chow rings}
Recall that $N_{L}=\overline{N}_{L}\oplus N_{L,tor}$, where $N_{L,tor}$ is the torsion subgroup
of $N_{L}$.  Let $N_{\Sigma_{\theta}}=N_{L,tor}\oplus |\Sigma_{\theta}|$.  By Theorem 1.1, we have 

\begin{prop}\Label{orbifoldlawrence}
The orbifold Chow ring $A^{*}_{orb}(\mathcal{X}(\mathbf{\Sigma}_\theta))$ of the 
Lawrence toric  DM stack $\mathcal{X}(\mathbf{\Sigma_{\theta}})$ is isomorphic to the ring
\begin{equation}\Label{chowringlawrence}
\frac{\mathbb{Q}[N_{\Sigma_{\theta}}]}
{\{\sum_{i=1}^{m}e(b_{L,i})y^{b_{L,i}}+\sum_{i=1}^{m}e(b'_{L,i})y^{b'_{L,i}}:e\in N_{L}^{\star}\}}.
\end{equation}~ $\square$
\end{prop}

Recall in \cite{JT} that  for any $c\in N$, there is a cone 
$\sigma\in \Delta_\mathbf{\beta}$ such that 
$\overline{c}=\sum_{\rho_{i}\subseteq \sigma}\alpha_{i}\overline{b}_{i}$ where 
$\alpha_{i}>0$ are  rational numbers. Let $N^{\Delta_\mathbf{\beta}}$ denote  
all the pairs $(c,\sigma)$. Then $N^{\Delta_\mathbf{\beta}}$ gives rise a 
group ring
$$\mathbb{Q}[\Delta_\mathbf{\beta}]=\bigoplus_{(c,\sigma)\in N^{\Delta_\mathbf{\beta}}}\mathbb{Q}\cdot y^{(c,\sigma)},$$
where $y$ is a formal variable.  
For any $(c,\sigma)\in N^{\Delta_\mathbf{\beta}}$, there exists a unique element
$(v,\tau)\in Box(\Delta_\mathbf{\beta})$ such that $\tau\subseteq\sigma$ and 
$c=v+\sum_{\rho_{i}\subseteq \sigma}m_{i}b_{i}$, 
where 
$m_{i}$ are nonnegative integers.  We call $(v,\tau)$ the fractional part of $(v,\sigma)$. We define the $ceiling ~function$ for fans. 
For $(c,\sigma)$ define the ceiling function $\lceil c \rceil_{\sigma}$ by  
$\lceil c \rceil_{\sigma}=\sum_{\rho_{i}\subseteq \tau}b_{i}+\sum_{\rho_{i}\subseteq \sigma}m_{i}b_{i}$. Note that 
if $\overline{v}=0$, $\lceil c \rceil_{\sigma}=\sum_{\rho_{i}\subseteq \sigma}m_{i}b_{i}$.  
For two pairs  $(c_1,\sigma_1)$, $(c_2,\sigma_2)$, if $\sigma_1\cup\sigma_2$ is a cone in $\Delta_\mathbf{\beta}$, define 
$\epsilon(c_1,c_2):=\lceil c_1 \rceil_{\sigma_{1}}+\lceil c_2 \rceil_{\sigma_{2}}-\lceil c_1+c_2 \rceil_{\sigma_{1}\cup\sigma_2}$.
Let $\sigma_{\epsilon}\subseteq\sigma_1\cup\sigma_2$ be the minimal cone in $\Delta_\mathbf{\beta}$ containing $\epsilon(c_1,c_2)$ so that 
$(\epsilon(c_1,c_2),\sigma_{\epsilon})\in N^{\Delta_\mathbf{\beta}}$. We define the grading on $\mathbb{Q}[\Delta_{\mathbf{\beta}}]$ as follows.
For any  $(c,\sigma)$, write $c=v+\sum_{\rho_{i}\subseteq \sigma}m_{i}b_{i}$,  then
$deg(y^{(c,\sigma)})=|\tau|+\sum_{\rho_{i}\subseteq\sigma}m_{i}$, where  
$|\tau|$ is the dimension of $\tau$.
By abuse of notation, we write $y^{(b_{i},\rho_i)}$ as $y^{b_{i}}$.
The multiplication 
is defined by 
\begin{equation}\Label{product1}
y^{(c_{1},\sigma_{1})}\cdot y^{(c_{2},\sigma_{2})}:=
\begin{cases}
(-1)^{|\sigma_{\epsilon}|}y^{(c_{1}+c_{2}+\epsilon(c_1,c_2),\sigma_1\cup\sigma_2)}&\text{if
$\sigma_{1}\cup\sigma_{2}$ is a cone in $\Delta_{\mathbf{\beta}}$}\,,\\
0&\text{otherwise}\,.
\end{cases} 
\end{equation}
From the property of $ceiling~ function$ we check that the multiplication is commutative and associative. So $\mathbb{Q}[\Delta_\mathbf{\beta}]$ is 
a unital associative commutative	ring. 
In \cite{JT}, it is shown that 
\begin{equation}\Label{chowringhyper}
A^{*}_{orb}(\mathcal{M}(\mathcal{A}))\cong \frac{\mathbb{Q}[\Delta_{\mathbf{\beta}}]}{\{\sum_{i=1}^{m}e(b_{i})y^{b_{i}}: e\in N^{\star}\}}.
\end{equation}

Consider the map $\beta: \mathbb{Z}^{m}\rightarrow N$ which is given by the vectors
$\{b_{1},\cdots,b_{m}\}$.  
We take $\{1,\cdots,m\}$ as the vertex set of the {\em matroid complex}
$M_{\beta}$, defined from $\beta$ by requiring that $F\in
M_{\beta}$ iff the  vectors
$\{\overline{b}_{i}\}_{i\in F}$ are linearly independent in $\overline{N}$. 
A face $F\in M_{\beta}$ corresponds 
to a cone in $\Delta_{\mathbf{\beta}}$ generated by $\{\overline{b}_{i}\}_{i\in F}$.  By
\cite{S}, the ``Stanley-Reisner'' ring of the matroid
$M_{\beta}$ is
$$\mathbb{Q}[M_{\beta}]=\frac{\mathbb{Q}[y^{b_{1}},\cdots,y^{b_{m}}]}{I_{M_{\beta}}},$$
where $I_{M_{\beta}}$ is the matroid ideal generated by the set
of  square-free monomials
$$\{y^{b_{i_{1}}}\cdots y^{b_{i_{k}}}|
\overline{b}_{i_{1}},\cdots,\overline{b}_{i_{k}} ~\text{linearly
dependent in}~\overline{N}\}.$$  
It is proved in  \cite{JT} that ,
$$\mathbb{Q}[\Delta_{\mathbf{\beta}}]\cong\bigoplus_{(v,\sigma)\in Box(\Delta_{\mathbf{\beta}})}y^{(v,\sigma)}\cdot \mathbb{Q}[M_{\beta}].$$
For any  $(v_{1},\sigma_{1}),(v_{2},\sigma_{2})\in
Box(\Delta_{\mathbf{\beta}})$, 
let $(v_{3},\sigma_{3})$ be the unique element in $Box(\Delta_{\mathbf{\beta}})$
such that $v_1+v_2+v_3\equiv 0$ in the local group given by $\sigma_1\cup\sigma_2$, where 
$\equiv 0$ means that there exists a cone  $\sigma(\overline{v}_{1},
\overline{v}_{2},\overline{v}_{3})$ in
$\Delta_{\mathbf{\beta}}$ such that 
$\overline{v}_{1}+\overline{v}_{2}+\overline{v}_{3}=\sum_{\rho_{i}\subseteq\sigma(\overline{v}_{1},
\overline{v}_{2},\overline{v}_{3})}a_{i}\overline{b}_{i}$,
where  $a_{i}=1 ~\text{or}~ 2$.  Let $\overline{v}_1=\sum_{\rho_j\subseteq\sigma_1}\alpha_{j}^{1}\overline{b}_{j}$, 
$\overline{v}_2=\sum_{\rho_j\subseteq\sigma_2}\alpha_{j}^{2}\overline{b}_{j}$,
$\overline{v}_3=\sum_{\rho_j\subseteq\sigma_3}\alpha_{j}^{3}\overline{b}_{j}$ with
$0<\alpha_{j}^{1},\alpha_{j}^{2},\alpha_{j}^{3}<1$. Let $I$ be the set of $i$ such that
$a_{i}=1$ and $\alpha_{j}^{1},\alpha_{j}^{2},\alpha_{j}^{3}$ exist, $J$  the set of $j$ such that $\rho_{j}$ belongs
to $\sigma(\overline{v}_{1},
\overline{v}_{2},\overline{v}_{3})$ but not $\sigma_{3}$.
If $(v,\sigma)\in Box(\Delta_\mathbf{\beta})$, let 
$(\check{v},\sigma)$ be the inverse of $(v,\sigma)$. Except torsion elements, equivalently, if $\overline{v}=\sum_{\rho_{i}\subseteq \sigma}\alpha_{i}\overline{b}_{i}$ for
$0<\alpha_{i}<1$, then $\check{\overline{v}}=\sum_{\rho_{i}\subseteq \sigma}(1-\alpha_{i})\overline{b}_{i}$. 
By abuse of notation, we write $y^{(b_{i},\rho_i)}$ as $y^{b_{i}}$.
We have that  $v_{1}+v_{2}=\check{v}_{3}+\sum_{a_{i}=2}b_{i}+\sum_{j\in J}b_{j}$. From (\ref{product1}), Lemma 5.11
and Lemma 5.12 in \cite{JT}, 
if $\overline{v}_{1},\overline{v}_{2}\neq 0$, we have 
$$
\lceil v_1\rceil_{\sigma_1}+\lceil v_2\rceil_{\sigma_2}-\lceil v_1+v_2\rceil_{\sigma_1\cup\sigma_2}=
\begin{cases}
\sum_{i\in I}b_{i}+\sum_{j\in J}b_{j}&\text{if
$\overline{v}_{1}\neq\check{\overline{v}}_{2}$}\,,\\
\sum_{j\in J}b_{j}&\text{if
$\overline{v}_{1}=\check{\overline{v}}_{2}$}\,.\\
\end{cases}
$$ 
So it is easy to check that the multiplication 
$y^{(v_{1},\sigma_{1})}\cdot y^{(v_{2},\sigma_{2})}$ can be written as
\begin{equation}\Label{product}
\begin{cases}
(-1)^{|I|+|J|}y^{(\check{v}_{3},\sigma_{3})}\cdot\prod_{a_{i}=2}
y^{b_{i}}\cdot\prod_{i\in I}
y^{b_{i}}\cdot \prod_{j\in J}y^{2b_{j}}&\text{if
$\overline{v}_{1},\overline{v}_{2}\in\sigma$ for $\sigma \in\Delta_{\mathbf{\beta}}$ and $\overline{v}_{1}\neq \check{\overline{v}}_{2}$}\,,\\
(-1)^{|J|} \prod_{j\in J}y^{2b_{j}}&\text{if
$\overline{v}_{1},\overline{v}_{2}\in\sigma$ for $\sigma \in\Delta_{\mathbf{\beta}}$ and $\overline{v}_{1}=\check{\overline{v}}_{2}$}\,,\\
0&\text{otherwise}\,.
\end{cases} 
\end{equation}

The following is the main result of this Section.

\begin{thm}\Label{main2}
There is an isomorphism of orbifold Chow rings $A_{orb}^{*}(\mathcal{X}(\mathbf{\Sigma_{\theta}}))\cong A_{orb}^{*}(\mathcal{M}(\mathcal{A}))$.
\end{thm}

\begin{pf}
The ring $\mathbb{Q}[N_{\Sigma_{\theta}}]$ is generated by $\{y^{b_{L,i}},y^{b'_{L,i}}: i=1,\cdots,m\}$ 
and $y^{v_{\theta}}$ for $v_{\theta}\in Box(\mathbf{\Sigma_{\theta}})$ by the definition. By Lemma \ref{box}, define a morphism
$$\phi: \mathbb{Q}[N_{\Sigma_{\theta}}]\rightarrow\mathbb{Q}[\Delta_{\mathbf{\beta}}] $$
by $y^{b_{L,i}}\mapsto y^{b_{i}}, y^{b'_{L,i}}\mapsto -y^{b_{i}}$ and $y^{v_{\theta}}\mapsto y^{(v,\sigma)}$. By \cite{HS}, the ideal $\mathcal{I}_{\theta}$ goes to the ideal $I_{M_{\beta}}$ and the relation 
$\{\sum_{i=1}^{m}e(b_{L,i})y^{b_{L,i}}+\sum_{i=1}^{m}e(b'_{L,i})y^{b'_{L,i}}:e\in N_{L}^{\star}\}$
goes to the relation $\{\sum_{i=1}^{m}e(b_{i})y^{b_{i}}:e\in N^{\star}\}$. Thus the two rings are isomorphic
as modules. 

It remains to check the multiplications. For any $y^{v_{\theta}}$ and $y^{b_{L,i}}$ or $y^{b'_{L,i}}$, let 
$y^{(v,\sigma)}$ be the corresponding element in $\mathbb{Q}[\Delta_{\mathbf{\beta}}]$. 
By the property of $v_{\theta}$ and Lemma \ref{box}, the minimal cone in 
$\Sigma_{\theta}$ containing $\overline{v}_{\theta}, \overline{b}_{L,i}$ must contains $\overline{b}'_{L,i}$.
By Lemma \ref{conemulti}, there is a cone in $\Delta_{\mathbf{\beta}}$ containing $\overline{v}, \overline{b}_{i}$. 
In this way, $y^{v_{\theta}}\cdot y^{b_{L,i}}$ goes to $y^{(v,\sigma)}\cdot y^{b_{i}}$ and 
$y^{v_{\theta}}\cdot y^{b'_{L,i}}$ goes to $-y^{(v,\sigma)}\cdot y^{b_{i}}$.
If there is no cone in $\Sigma_{\theta}$
containing $\overline{v}_{\theta}, \overline{b}_{L,i},\overline{b}'_{L,i}$, then by Lemma \ref{conemulti} there is no cone in $\Delta_{\mathbf{\beta}}$ containing $\overline{v}, \overline{b}_{i}$. So 
$y^{v_{\theta}}\cdot y^{b_{L,i}}=0$ goes to $y^{(v,\sigma)}\cdot y^{b_{i}}=0$ and 
$y^{v_{\theta}}\cdot y^{b'_{L,i}}=0$ goes to $-y^{(v,\sigma)}\cdot y^{b_{i}}=0$.

For any $y^{v_{\theta,1}},y^{v_{\theta,2}}$, let $y^{(v_{1},\sigma_{1})},y^{(v_{2},\sigma_{2})}$ be the corresponding 
elements in $\mathbb{Q}[\Delta_{\mathbf{\beta}}]$. If there is no cone 
in $\Sigma_{\theta}$ containing $\overline{v}_{\theta,1},\overline{v}_{\theta,2}$, then by Lemmas \ref{conemulti}
and \ref{box}, there is no cone in $\Delta_{\mathbf{\beta}}$ containing $\overline{v}_{1},\overline{v}_{2}$. So
$y^{v_{\theta,1}}\cdot y^{v_{\theta,2}}=0$ goes to $y^{(v_{1},\sigma_{1})}\cdot y^{(v_{2},\sigma_{2})}=0$.
Suppose there
is a cone  containing $\overline{v}_{\theta,1},\overline{v}_{\theta,2}$, let 
$v_{\theta,3}\in Box(\mathbf{\Sigma}_{\theta})$ such that $v_{\theta,1}+v_{\theta,2}+v_{\theta,3}\equiv 0$.
Let $\sigma(\overline{v}_{\theta,1},
\overline{v}_{\theta,2},\overline{v}_{\theta,3})$ be the minimal cone containing $\overline{v}_{\theta,1},
\overline{v}_{\theta,2},\overline{v}_{\theta,3}$ in $\Sigma_{\theta}$. Then by Lemmas \ref{conemulti}
and \ref{box}, 
$\sigma(\overline{v}_{\theta,1},
\overline{v}_{\theta,2},\overline{v}_{\theta,3})$ is the Lawrence lifting of $\sigma(\overline{v}_{1},
\overline{v}_{2},\overline{v}_{3})$ for  $(v_{1},\sigma_{1}),(v_{2},\sigma_{2}),(v_{3},\sigma_{3})\in Box(\Delta_{\mathbf{\beta}})$.
So we may write
$\overline{v}_{\theta,1}+\overline{v}_{\theta,2}+\overline{v}_{\theta,3}=\sum_{\rho_{i}\subseteq \sigma(\overline{v}_{1},\overline{v}_{2},\overline{v}_{3})}a_{i}\overline{b}_{L,i}
+\sum_{\rho_{i}\subseteq \sigma(\overline{v}_{1},\overline{v}_{2},\overline{v}_{3})}a'_{i}\overline{b}'_{L,i}$.
The corresponding $\overline{v}_{1}+\overline{v}_{2}+\overline{v}_{3}=
\sum_{\rho_{i}\subseteq \sigma(\overline{v}_{1},\overline{v}_{2},\overline{v}_{3})}a_{i}\overline{b}_{i}$. 
Let
$(\check{v},\sigma)$ be  the inverse of $(v,\sigma)$ in
$Box(\Delta_{\mathbf{\beta}})$, i.e. if $v$ is nontorsion and
$\overline{v}=\sum_{\rho_{i}\subseteq \sigma}\alpha_{i}\overline{b}_{i}$
for $0<\alpha_{i}<1$, then  $\check{\overline{v}}=\sum_{\rho_{i}\subseteq
\sigma}(1-\alpha_{i})\overline{b}_{i}$. The $\check{v}_{\theta}$ is defined similarly
in $Box(\mathbf{\Sigma_{\theta}})$. The notation  $J$ represents the set of $j$ such that $\rho_{j}$ belongs
to $\sigma(\overline{v}_{1},\overline{v}_{2},\overline{v}_{3})$ but not $\sigma_{3}$,   the 
corresponding $\rho_{L,j},\rho'_{L,j}$ belong
to $\sigma(\overline{v}_{\theta,1},\overline{v}_{\theta,2},\overline{v}_{\theta,3})$ but not $\sigma(\overline{v}_{\theta,3})$.

If some $\overline{v}_{\theta,i}=0$ which means that $v_{\theta,i}$ is a torsion. Then 
from Lemma (\ref{box}) the corresponding 
$v$ is also a torsion element. In this case we know that the orbifold cup product 
$y^{v_{\theta,1}}\cdot y^{v_{\theta,2}}$ is the usual product, and under the 
map $\phi$,  is equal to
$y^{(v_{1},\sigma_{1})}\cdot y^{(v_{2},\sigma_{2})}$.

If $\overline{v}_{\theta,1}=\check{\overline{v}}_{\theta,2}$, then $\overline{v}_{\theta,3}=0$ and the obstruction bundle over the corresponding 3-twisted sector is zero. 
The set $J$ is the set $j$ such that $\rho_j$ belongs to $\sigma(\overline{v}_{\theta,1})$.  
So from \cite{BCS}, we have 
$$y^{v_{\theta,1}}\cdot y^{v_{\theta,2}}=
\prod_{j\in J}y^{b_{L,j}}\cdot y^{b'_{L,j}}.$$
Under the map $\phi$ we see that $y^{(v_{1},\sigma_{1})}\cdot y^{(v_{2},\sigma_{2})}$ is equal to the second line in the  product 
(\ref{product}). 

If $\overline{v}_{\theta,1}\neq\check{\overline{v}}_{\theta,2}$, then $\overline{v}_{\theta,3}\neq 0$ and the obstruction bundle
is given by Proposition \ref{obstructionbdle}.  If all  $\alpha_{j}^{1},\alpha_{j}^{2},\alpha_{j}^{3}$ exist, the coefficients $a_{i}$ and $a'_{i}$ satisfy that if $a_{i}=1$ then
$a'_{i}=2$, and if $a_{i}=2$ then
$a'_{i}=1$. So
from \cite{BCS},
$$y^{v_{\theta,1}}\cdot y^{v_{\theta,2}}=
y^{\check{v}_{\theta,3}}\cdot\prod_{a_{i}=2}
y^{b_{L,i}}\cdot\prod_{i\in I}
y^{b'_{L,i}}\cdot \prod_{j\in J}y^{b_{L,j}}\cdot y^{b'_{L,j}}.$$
Under the map $\phi$ we see that $y^{(v_{1},\sigma_{1})}\cdot y^{(v_{2},\sigma_{2})}$ is equal to the first line in the  product 
(\ref{product}). By Lemma \ref{box}, the box elements have the same orbifold degrees. 
By Corollary \ref{sectors} and the definition of orbifold cup product in (\ref{cupproduct}),
the products $y^{v_{\theta,1}}\cdot y^{v_{\theta,2}}$ and $y^{(v_{1},\sigma_{1})}\cdot y^{(v_{2},\sigma_{2})}$ have 
the same degrees in both Chow rings.
So $\phi$ induces a ring isomorphism 
$A_{orb}^{*}(\mathcal{X}(\mathbf{\Sigma_{\theta}}))\cong A_{orb}^{*}(\mathcal{M}(\mathcal{A}))$.
\end{pf}

\begin{rmk}
The presentation (\ref{chowringhyper}) of orbifold Chow ring only depends on the matroid complex corresponding to the map 
$\beta: \mathbb{Z}^{m}\rightarrow N$, not $\theta$. Note that the presentation (\ref{chowringlawrence}) depends on 
the fan $\Sigma_{\theta}$. We couldn't see explicitly from this presentation that the ring is independent to the 
choice of generic elements $\theta$.
\end{rmk}


\end{document}